\documentclass[12pt]{article}
\usepackage{amsmath,amsthm,amssymb,latexsym,color}
\usepackage{bbm}
\title{Random polytopes obtained by matrices with heavy tailed entries.}
\author{O. Gu\'edon  \and A. E. Litvak \and K. Tatarko
}

\date{ }

\newcommand\address{\noindent\leavevmode

\medskip

\noindent
Olivier Gu\'{e}don \\
Universit\'{e} Paris-Est Marne-La-Vall\'ee\\
Laboratoire d'Analyse et de Math\'{e}matiques Appliqu\'ees, \\
5, boulevard Descartes,
Champs sur Marne,\\
77454 Marne-la-Vall\'{e}e,  Cedex 2, France\\
\texttt{\small%
e-mail:  olivier.guedon@univ-mlv.fr}\\

\medskip

\noindent
A. E. Litvak
and K. Tatarko,\\
Dept.~of Math.~and Stat.~Sciences,\\
University of Alberta, \\
Edmonton, AB, Canada, T6G 2G1.\\
\texttt{\small
e-mails:  aelitvak@gmail.com \, \, and \, \,
ktatarko@gmail.com}
}

\def\no{ \nor \cdot  \nor }

\def\nx{ \nor x \nor }

\def\ny{ \nor y \nor }

\def\kkk{\null\hfill $\Box $ \\}
\def\r{ \right}

\def\eps{\varepsilon}

\def\lam{\lambda}

\def\nor{ \| }

\def\gam{\gamma}

\newcommand{\RR}{\mathbb{R}}

\newcommand{\R}{\RR ^n}

\newcommand{\pp}{\mathbb{P}}
\newcommand{\E}{\mathbb{E}}
\newcommand{\EE}{\mathbb{E}}

\newcommand{\K}{b}
\newcommand{\KK}{M}

\newcommand{\proofs}{\noindent {\bf Proof. }}

\newcommand{\nk}[1]{\| #1 \|_{k,2}}
\newcommand{\Bnk}{\B_{k,2}}
\newcommand{\Xnk}{X_{k,2}}

\newcommand{\euclidnorm}[1]{\| #1 \|_2} 

\def\la{\left\langle}
\def\ra{\r\rangle}

\newtheorem{theorem}{Theorem}[section]
\newtheorem{lemma}[theorem]{Lemma}
\newtheorem{remark}[theorem]{Remark}

\newtheorem{prop}[theorem]{Proposition}
\newtheorem{sled}[theorem]{Corollary}

\newcommand{\conv}{\mathop{\rm conv\,}}
\newcommand{\absconv}{\mathop{\rm abs \,conv\,}}

\renewcommand{\P}{\mathbb{P}}

\newcommand{\B}{{\textbf B}}

\def\cf{\mathcal{Q}}

\begin{document}

\maketitle

\begin{abstract}
 Let $\Gamma$ be an $N\times n$ random matrix with
 independent entries and such that
in each row entries are i.i.d. Assume also that the entries are
symmetric, have unit
variances, and satisfy  a small ball probabilistic estimate uniformly.
We investigate properties of the corresponding random  polytope
$\Gamma^* B_1^N$ in $\R$ (the absolute convex hull of rows of $\Gamma$).
In particular, we show that
$$
 \Gamma^* B_1^N \supset \K^{-1} \left( B_{\infty}^n \cap \sqrt{\ln (N/n)}\,  B_2^n \r),
$$
where $\K$ depends only on parameters in small ball inequality.
 This extends results of \cite{LPRT} and recent results of \cite{KKR}.
This inclusion is equivalent to so-called $\ell_1$-quotient property
and plays an important role in compressed  sensing
(see \cite{KKR} and references therein).
\end{abstract}

\smallskip

\noindent
{\small \bf AMS 2010 Classification:}
{\small
primary: 52A22, 46B06, 60B20,
secondary: 52A23, 46B09, 15B52.
}


\section{Introduction}
In this note, we deal with a rectangular $N \times n$ random matrices
$\Gamma = \{\xi_{ij}\}_{1 \le i \le N, 1 \le j \le n}$, where
$\xi_{ij}$ are independent symmetric random variables with unit variance
satisfying uniform small ball probabilistic estimate.
More precisely, in the main theorem we assume that there exist
$u, v \in (0,1)$ such that
\begin{equation}\label{smalbal}
 \forall i, j \quad \sup_{\lambda \in \RR} \pp \big\{ |\xi_{ij} - \lambda | \le u \big\} \le v.
\end{equation}
Of course, if variables have a bounded moment $r>2$, we will
have  better  estimates.
 We are interested in geometric parameters of the random polytope
generated by $\Gamma$, that is, the absolute convex hull of
rows of $\Gamma$. In other words, the random polytope under consideration is
$\Gamma^* B_1^N$, where $B_1^N$ is the $N$-dimensional octahedron (cross-polytope).
 Such random polytopes have been extensively
studied in the literature, especially in the Gaussian case and
in the Bernoulli case. The Gaussian random polytopes
in the case when $N$ is proportional to $n$ have many applications in
the Asymptotic Geometric Analysis (see e.g., \cite{gldistance}
and \cite{szbasis}, and the survey \cite{mat}).  The Bernoulli case
corresponds to 0/1 random polytopes. For their combinatorial properties
we refer the reader to  \cite{Furedi, Barany} (see also the survey \cite{Ziegler}).
Their  geometric parameters have been studied in \cite{GH, LPRT}.
In the compressed sensing it was shown that
the so-called $\ell_1$-quotient property is responsible for robustness
in  certain $\ell_1$-minimizations (see  \cite{KKR} and references
therein). More precisely, an $n\times N$ (with $N\geq n$) matrix $A$ satisfies
the $\ell_1$-quotient property with a constant $\K$ relative to
a norm $\no$ if for every $y\in \RR^n$ there exists $x\in \RR^N$
such that $Ax = y$ and $\nx_1\leq \K \sqrt{n/\ln (e N/n)}\, \ny$, where
$\no_1$ denotes the $\ell_1$-norm. It is easy to see that geometrically this
means
$$
      B_{\no} \subset \K\sqrt{n/\ln (e N/n)}\, A B_1^N,
$$
where $B_{\no}$ is the unit ball of $\no$. To prove robustness of
noise-blind compressed sensing, the authors of \cite{KKR} dealt with
the norm
$$
 \no = \max\{\|\cdot\|_2, \sqrt{\ln (e N/n)}\no_\infty\},
$$
where
$\|\cdot\|_2$ is the standard Euclidean norm and $\no_\infty$ is the  $\ell_\infty$-norm.
Theorem~5 in \cite{KKR} states that assuming
that entries of $A$ are symmetric i.i.d. random variables with unit variances,
and that they have regular (in fact,  $\psi_\alpha$) behaviour of all moments
till the moment of order $\ln n$, the matrix $A/\sqrt{n}$  has the $\ell_1$-quotient property
with high probability.
Geometrically this means
\begin{equation}\label{incint}
   A B_1^N \supset \K^{-1} \left( B_{\infty}^n \cap \sqrt{\ln (N/n)}\,  B_2^n \r).
\end{equation}
The work \cite{KKR} complements results of \cite{LPRT}, where this inclusion
was proved for random matrices with symmetric i.i.d. entries having at
least third bounded moment and such that the operator
norm of the matrix is bounded  with high probability.

The main purpose of this note is to prove such an inclusion with
much weaker assumptions on the distribution of the entries. In fact,
we require only boundedness of  second moments. Thus ``robustness" Theorem~8 in \cite{KKR}  holds
under much weaker assumptions on the random matrix. Our main result is the following
theorem (see Theorem~\ref{inclusion} for slightly better probability estimates).

 \begin{theorem}\label{main-int}
 There exist positive constants $\K, \KK$ depending only on  $u$ and $v$
 and an absolute  constant $c >0$ such that the following holds.
 Let $N\geq \KK n$ and  assume that the entries of an $N\times n$ random matrix $\Gamma$ are
   independent symmetric random variables with unit variances satisfying condition
   (\ref{smalbal}) and such that in each row the entries are i.i.d.
   Then with probability at least $1-\exp(-c n)$ the inclusion (\ref{incint})
   holds for the matrix  $A=\Gamma^*$.
 \end{theorem}

We use this theorem to study geometric properties of random polytopes
$K_N=\Gamma^*B_1^N$ and $K_N^0$, such as behavior of their volumes and mean widths.
Our ``volume" theorem states the following (see Theorems~\ref{volone} and
\ref{voltwo} for more precise statements).

 \begin{theorem}\label{vol-int}
 There exist positive constants $C_1, C_2$ depending only on  $u$ and $v$
 and  absolute positive constants $C,c$ such that for $C_1 n\leq N\leq e^n$
  with probability at least $1-\exp(-c n)$ one has
 $$
|K_N|^{1/n} \geq  C_2
       \sqrt{\frac{ \ln (N/n)}{n}} \, \, \, \,\, \,\, \,
\mbox{ and } \, \, \, \,\, \, \, \, |K_N^0|^{1/n} \leq
\frac{C}{C_2  \sqrt{  n \ln (N/n)}},
$$
 where $K_N=\Gamma^*B_1^N$ and the  matrix $\Gamma$ is as in Theorem~\ref{main-int}.
 Moreover, the bounds on the volumes are sharp, provided that the Euclidean
 lengths of the rows of $\Gamma $ are of order of $\sqrt{n}$ at most.
 \end{theorem}

Our proof of Theorem~\ref{main-int} follows the general scheme of \cite{LPRT} with a very delicate change --
in \cite{LPRT} there was an assumption that the operator
norm
of $\Gamma$ is bounded by $C\sqrt{N}$ with high probability. However it is known
that such a bound does not hold in general unless fourth moments are bounded
(\cite{silv}, see also \cite{LiSp} for quantitative bounds).
To avoid using the norm of $\Gamma$, we use ideas appearing
in \cite{RT}, where the authors constructed a certain deterministic $\eps$-net
(in $\ell_2$-metric) $\mathcal{N}$ such that $A \mathcal{N}$ is a good net for $A B_2^n$
for most realizations of a square random  matrix $A$. We extend their construction
in three directions. First, we work with rectangular random matrices, not only square matrices.
Second, we need a net for the image of a given convex body (not only for the image of the unit
Euclidean ball). Finally, instead of approximation in the Euclidean norm only, we use
approximation in the following norm
\begin{equation}
\label{nk}
\nk{a} = \left(\sum_{i=1}^k (a_i^*)^2 \right)^{1/2},
\end{equation}
where $1\leq k\leq N$ and $a_1^* \ge a_2^* \ge \ldots \ge a_N^*$
is the decreasing rearrangement
of the sequence of numbers $|a_1|, \ldots, |a_N|$.
This norm appears naturally and plays a crucial role in our proof of inclusion (\ref{incint}).
The generalization of the net from \cite{RT} is a new key ingredient, see Theorem~\ref{Kostya}.
We would like to emphasize, that norms $\nk{\cdot}$ played an important role in proofs of
many results of Asymptotic Geometric Analysis, see e.g. \cite{Gl3, GGMP, OG}. For the systematic
studies of norms  $\nk{\cdot}$ and their unit balls we refer to \cite{GLSW}. We believe
that the new approximation in $\nk{\cdot}$ norms will  find other applications in the theory.
In the last section, we present one more application of
Theorem~\ref{Kostya} -- we show that it can be used to estimate the  smallest singular value
of a tall random matrix -- see the discussion at the beginning of Section~\ref{smsing}.

\bigskip

\noindent
{\bf Acknowledgement. }
This project has been started when the second named author
was visiting University Paris-Est at Marne-la-Vall\'ee.
He is thankful for excellent working conditions.
All three authors are grateful to MFO, Oberwolfach, where
part of the work was done during the workshop ``Convex Geometry and its Applications".

\section{Notations}
\label{prem}

By $\la \cdot , \cdot \ra$ we denote the canonical inner product on the $m$-dimensional real space $\RR ^m$
and for $1\leq p \leq \infty$,   the $\ell _p$-norm is defined for any $a \in \RR^m$ by
$$
\|a\| _p = \left( \sum _{i = 1}^m |a_i|^p \r) ^{1/p} \, \mbox{ for
 } \ p < \infty \, \, \, \,  \, \, \, \, \mbox{ and } \, \, \, \,  \, \, \, \, \|a\| _{\infty }
 = \sup _{i = 1, \ldots, m} |a_i|.
$$
As usual, $\ell _p^m = (\RR ^m, \|\cdot \|_p)$, and the unit ball of
$\ell _p^m$ is denoted by $B_p^m$. The unit sphere of $\ell _2^m$ is
denoted by $S^{m-1}$, and   the canonical basis of $\ell _2^m$
is  denoted  by  $e_1, \ldots, e_m$.

Given an integer $k \in \{1, \ldots, N\}$, we denote by $\Xnk$ the normed space  $\RR^N$
equipped with the norm $\nk{\cdot}$ defined by (\ref{nk}).
The unit ball of $\Xnk$ is denoted by $\Bnk$. Note that for
$k = N$ we have $\nk{a} = \|a\|_2$ and that for any $k \le N$ and any  $a \in \RR^N$,
\[
\nk{a} \le \|a\|_2 \le \sqrt{\frac{N}{k}}\, \nk{a} \mathrm{ \quad \quad or, equivalentely,  \quad\quad}
B_2^N \subset \Bnk \subset \sqrt{\frac{N}{k}}\, B_2^N.
\]

Given  integers $\ell \geq k\geq 1$, we
denote $[k]=\{1, 2, \ldots, k\}$ and $[k, \ell]=\{k, k+1, \ldots, \ell\}$.
Given a number $a$ we denote the largest integer not exceeding $a$ by
$\lfloor a\rfloor $ and the smallest integer larger than or equal to $a$ by
$\lceil a \rceil$.

Given points $x_1, \ldots, x_k$ in $\RR ^m$ we denote their convex
hull by $\conv \{x_i\}_{i\leq k}$ and their absolute convex
hull by $\absconv \{x_i\}_{i\leq k} =
\conv \{\pm x_i\}_{i\leq k}$.
Given $\sigma \subset [m]$ by $P_{\sigma}$ we denote
the coordinate projection onto
$\RR ^{\sigma}=\{x\in \RR^m \, | \, x_i=0 \, \mbox{ for } \, i\notin \sigma\}$.

Given a finite set $E$ we denote its cardinality by $|E|$.
We also use
$|K|$ for the volume of a body $K\subset \RR^m$ (and, more generally,
for the $m$-dimensional Lebesgue measure of a measurable subset in
$\RR^m$).
Let $K$ be a symmetric convex body with non empty interior. We denote
its Minkowski's functional by $\| x \| _K$. The support function of $K$ is
$ h_K (x) = \sup _{y\in K} \la x, y\ra,$
the polar of $K$ is
$$
 K^0 = \left\{ x\in \RR^m \ | \ \la x, y \ra \leq 1\, \mbox{ for every }
 \, y \in K \r\}.
$$
Note that $h_K (\cdot) = \| \cdot \| _{K^0}$.

Given a set $L\subset \RR ^m$, a convex body $K\subset \RR^m$, and
$\eps >0$ we say that a subset $\mathcal{N} \subset \RR^m$
is an {\it $\eps$-net} of $L$ with
respect to $K$ if
$$
 \mathcal{N}\subset L \subset \bigcup _{x\in \mathcal{N}}\ (x + \eps K).
$$
The cardinality of the smallest $\eps$-net  of $L$ with
respect to $K$ we denote by $N(L, \eps K)$.

For a given probability space, we denote by $\P (\cdot)$ and $\E$ the probability
of an event and the expectation respectively.  A $\pm 1$ random variable
taking values $1$ and $-1$ with probability $1/2$ is called
a Rademacher random variable.

In this paper we are interested in rectangular $N\times n$ matrices
$\Gamma =\{\xi_{ij}\}_{\stackrel{1 \le i \le N}{1 \le j \le n}}$, with $N\ge n$, where the entries are real-valued random
variables on some probability space $(\Omega,{\mathcal{A}},\P)$.
We will mainly consider the following model of matrix  $\Gamma$:
\begin{equation}
\label{H}
\left\{
\begin{array}{l}
 \forall i,j \quad \xi_{ij} \  \mathrm{ are  \ independent, \ symmetric \ and } \  \E \xi_{ij}^2 = 1,
\\
 \mathrm{in \ each \ row \ the \  entries \ are \ identically \ distributed.}
 \end{array}
 \right.
\end{equation}
At the beginning of Section~\ref{geom}, we will also assume that the entries of $\Gamma$ satisfy a uniform small ball estimate.
If  $\xi_{ij} \sim {\cal N}(0, 1)$ are independent Gaussian random variables we say that
$\Gamma$ is {\it a Gaussian random matrix}.

\section{Construction of a good deterministic net.}
In this section we present a key result of this paper.
Let $T$ be a subset  of $\R$, we aim at constructing a deterministic
net such that for every general random operator $\Gamma : \RR^n \to \RR^N$, with overwhelming  probability,
the image of the net by the random operator $\Gamma$ is a good approximation of $\Gamma T$. We show that we can quantify this approximation by  almost any norm $\nk{\cdot}$ defined in $\eqref{nk}$.
For integers $1 \leq n \leq N$ and for $0\leq \delta \leq 1$, set
\begin{equation}
\label{eq:F}
 F(\delta, n, N)  =
  \left\{
\begin{array}{ll}
 \left(32\delta N/n\r)^n & \mbox{ if }\,\, \delta \geq n/(2N),
\\
  \left(e n/(\delta N)\r)^{4\delta N}& \mbox{ if }\,\,  \delta \leq  n/(2N).
\end{array}
\right.
\end{equation}

\begin{theorem}\label{Kostya}
 Let $n \in [N]$, $0\leq \delta \leq 1$, $0<\eps\leq  1$. Let $k \in [N]$ such that $k \ln (eN/k) \ge n$.
 Let  $T$ be a non-empty subset of $\R$ and denote $M:=N(T, \eps B_\infty^n)$.
 Then there exists a set $\mathcal{N}\subset T$ and  a
 collection of  parallelepipeds $\mathcal{P}$ in $\R$ such that
$$
 \max\{|\mathcal{N}|,  | \mathcal{P}|\} \leq
  M \, F(\delta, n, N)
   \, e^{\delta N}.
$$
Moreover, for any random matrix $\Gamma $ satisfying assumption \eqref{H},  with probability
at least $1- e^{- k \ln\left({eN}/{k}\right)}-e^{-\delta N/4}$, one has
$$
\left\{
\begin{array}{l}
 \forall x\in T \,\, \exists y\in \mathcal{N} \, \, \, \, \mbox{ such that } \, \, \, \,
\displaystyle
  \nk{\Gamma (x-y)} \le  C\eps  \sqrt{ \frac{k n}{\delta} \ln\left(\frac{eN}{k}\right)},
\\
\displaystyle
  \forall x\in T \,\, \exists P\in \mathcal{P} \, \, \, \, \mbox{ such that } \, \, \, \,
  x\in P \, \, \mbox{ and } \, \, \Gamma P \subset  \Gamma x +  C\eps \sqrt{ \frac{k n}{\delta}  \ln\left(\frac{eN}{k}\right)} \ \Bnk,
  \end{array}
\right.
$$
where $C\geq 1$ is an absolute  constant.
\end{theorem}

\begin{remark}{\rm
This result extends Theorem~A
and Corollary~A from \cite{RT},
 where the authors considered the case of square matrices,
 $T = S^{n-1}$ and $k=N$, which corresponds to the approximation of $\Gamma x$
 in the Euclidean norm.}
 \end{remark}
\subsection{Basic facts about covering numbers and operator norms of random matrices.}
We begin by recalling some classical estimates for covering numbers that will be used later.
It is well known that for any two centrally symmetric bodies $K$ and $L$ in $\RR^m$
and any $\eps >0$ there exists an $\eps$-net $\mathcal{N}$ of $L$ with respect to $K$
with cardinality
\begin{equation}\label{netcardgen}
 |\mathcal{N}|\leq |(2/\eps) L + K|/|K|
\end{equation}
  (see e.g. Lemma~4.16 in
\cite{P}).
In particular, if $K=L$ are centrally symmetric bodies in $\RR^m$
(or if $L$ is the boundary of a centrally symmetric body $K$) then
$|\mathcal{N}|\leq (1+2/\eps)^m$.
\begin{lemma}\label{covcubesch}
a) For every $\eps \in (0, 1/\sqrt{m}   ]$
$$
    N(B_2^m, \eps B_\infty^m) \leq \left(7/(\eps \sqrt{m})\r)^m
$$
and for every $\eps \in (1/\sqrt{m}, 1]$
$$
    N(B_2^m, \eps B_\infty^m) \leq \left(17 \eps^2 m \r)^{1/\eps ^2}.
$$
b) For $J \subset [m]$, let $S^J = \{x \in \RR^J\, | \,  \|x\|_2 = 1\}.$ For every $\eps \in (0,1)$ and every integer $k \le m$, there exists a finite set ${\cal N} \subset  \cup_{|J|=k} \, S^J$ such that
\begin{equation}
\label{eq:coversparse}
\begin{cases}
\displaystyle
| {\cal N} | \le
\exp{( k \ln(3/\eps) + k \ln(em/k))},
\\
\displaystyle
\forall J \subset [m] {\mathrm \ with  \ } |J|=k \,\,\,\,
\forall y \in S^J
\quad
\exists z \in {\cal N} \cap S^J  {\mathrm \ such  \ that \ }
\|y-z\|_2 \le \eps.
\end{cases}
\end{equation}
\end{lemma}

\proofs
a) Note that for every $m\geq 1$ one has  $(1/\sqrt{m} )  B_\infty^m \subset B_2^m$ and
$|B_2^m|\leq (2\pi e /m)^{m/2}$. Therefore,  by (\ref{netcardgen}), we obtain for every
$\eps  \leq 1/\sqrt{m}$
$$
  N(B_2^m, \eps B_\infty^m)\leq
  \frac{\left| \frac{2}{\eps} B_2^m + B_\infty^m \right|}{|B_\infty^m|}
  \leq \left(\frac{3}{\eps}\right)^m \frac{|B_2^m|}{|B_\infty^m|}
  \leq \left(\frac{3\sqrt{\pi e}}{\eps \sqrt{2m } }\r)^m.
$$
This implies the first bound. For the second bound
note that for every $x\in B_2^n$ the number of coordinates of $x$ larger than $\eps$
is at most $1/\eps^2$. Thus every $x\in B_2^n$ can be presented as $x=y+z$, where
the cardinality of support of $y$ is at most $1/\eps^2$, $z\in \eps B_\infty^n$,
and supports of $y$ and  $z$ are mutually disjoint. Therefore, it is enough to
cover $B_2^{\sigma}$ by $\eps B_\infty^{\sigma}$ for all $\sigma \subset [n]$ with
$|\sigma|=m := \lfloor 1/\eps^2\rfloor$. Using the above bound we obtain
$$
 N(B_2^n, \eps B_\infty^n) \leq \binom{n}{m}
 \left(\frac{3\sqrt{\pi e}}{\eps \sqrt{2m } }\r)^m
 \leq \left(\frac{3e n \sqrt{\pi e}}{\eps m \sqrt{2m } }\r)^m,
$$
which implies the desired result as $m\leq 1/\eps^2$.

b) Fix $\eps \in (0,1)$. For any fixed $J \subset  [m]$ of cardinality $k$, we cover $S^J$ by an $\eps$-net (of points in $S^J$) of cardinality at most $(1+2/\eps)^k \le (3/\eps)^k$ and we take the union of these nets over all sets $J$ of cardinality $k$. We conclude using that $\binom{m}{k} \le (em/k)^k$.
\kkk

The next lemma is a classical consequence of  estimates for  covering numbers for evaluating operator norms of random matrices.
\begin{lemma}\label{newlemma}
Let $B=\{b_{ij}\}_{\stackrel{1 \le i \le N}{1 \le j \le n}}$
be a fixed $N\times n$  matrix. Let $k \in [N]$ be
such that $k \ln \frac{eN}{k} \ge n$.
Let $\eps_{ij}$ be i.i.d. Rademacher random variables.
Denote $B_\eps = \{\eps_{ij} b_{ij}\}_{\stackrel{1 \le i \le N}{1 \le j \le n}}$.
Then for every $t\geq 1$ one has
$$
 \pp \left( \|B_\eps : \ell_\infty^n \to \Xnk \| \geq 6 \, t \,\sqrt{ k \ln\left(\frac{eN}{k}\right) }\,  \max_{i\leq N} \|R_i(B)\|_2 \r) \leq  e^{\displaystyle -t^2 k \ln\left(eN/k\right)},
$$
where $R_i(B)$, $i\leq N$, are the rows of $B$.
\end{lemma}

\proofs
Observe that for any $a \in \RR^N$, we have
\[
\nk{a} = \sup_{\stackrel{J \subset [N]}{|J|=k}} \ \sup_{b \in S^J} \sum_{i=1}^N a_i b_i.
\]
 Given $x\in \{\pm 1\}^n$, $y \in S^{N-1}$,  consider the following random variable,
$$
 \xi _{x,y} = \sum _{j=1}^{n} \sum_{i=1}^N \,\, \eps_{ij}  b_{ij}x_j y_i.
$$
Since $e^x+e^{-x}\leq 2 \exp(x^2/2)$ for every real $x$, we observe for $\lam >0$,
\begin{align*}
\E \exp\left(\lambda \sum _{j=1}^{n} \sum_{i=1}^N \,\, \eps_{ij}  b_{ij}x_j y_i\r)
& =
\prod_{j=1}^{n} \prod_{i=1}^N \E \exp \left( \lambda \eps_{ij}  b_{ij}x_j y_i\r)
\le \exp \left( \frac{\lambda^2}{2} \sum_{i=1}^N y_i^2 \|R_i(B)\|_2^2 \r)
\\
& \le \exp \left(\frac{\lambda^2}{2} \max_{i\leq N} \|R_i(B)\|_2^2\r).
\end{align*}
Therefore, using
the Laplace transform of $\xi_{x,y}$, we deduce that for any $u > 0$,
 \[
 \pp\left(  \xi _{x,y} > u  \max_{i\leq N} \|R_i(B)\|_2 \r) \le  e^{- u^2/2}.
 \]
Note that
\begin{equation}\label{matr-norm1}
 \|B_\eps : \ell_\infty^n \to \Xnk \|
   = \sup _{x\in \{\pm 1\}^n}  \sup_{\stackrel{J \subset [N]}{|J|=k}} \ \sup_{y \in S^J} \, \xi_{x,y}.
\end{equation}
Now we apply the classical net argument. Let $\cal N$ be the net defined by \eqref{eq:coversparse} with $\eps = 1/2$.
Then
\begin{align*}
\pp\left(\sup _{x\in \{\pm 1\}^n} \sup_{z \in {\cal N}}  \, \xi_{x,z}
\ge u \max_{i\leq N} \|R_i(B)\|_2 \r)
& \le  2^{n} |{\cal N}| e^{-u^2/2}
\\
& \le  2^{n} \exp\left( -\frac{u^2}{2} + k \ln 6 + k \ln(eN/k) \r).
\end{align*}
Taking $u = 3 t \sqrt{ k \ln(eN/k)}$ and using $k\ln (eN/k)\geq n$, we get for every $t \ge 1$,
\[
\pp\left(\sup _{x\in \{\pm 1\}^n} \sup_{z \in {\cal N}}  \, \xi_{x,z}
\ge 3 t \sqrt{ k \ln(eN/k)}\, \max_{i\leq N} \|R_i(B)\|_2 \r)
\le
 e^{\displaystyle -t^2 k \ln\left(eN/k\right)}.
\]
By definition of $\cal N$, for any $J \subset [N]$ of cardinality $k$ and $y \in S^J$, there exists $z \in {\cal N} \cap S^J$ such that $\|z-y\|_2 \le 1/2$, hence, by the triangle inequality,
\[
\sup _{x\in \{\pm 1\}^n}  \sup_{\stackrel{J \subset [N]}{|J|=k}} \ \sup_{y \in S^J} \, \xi_{x,y}
\le 2
\sup _{x\in \{\pm 1\}^n} \sup_{z \in {\cal N}}  \, \xi_{x,z}.
\]
This completes the proof of the lemma.
\kkk

\subsection{Auxiliary statements}
By  $\mathcal{D}_n$ we denote the set of all $n\times n$ diagonal matrices
whose diagonal entries belong to the set $\{1\}\cup \{ 2^{-2^{k}}\}_{k\geq 0}$.
The following theorem was proved in \cite{RT} in the square case.
However the proof works as well in the rectangular case.
One just needs to
repeat the proof of Proposition~2.7 there for $N\times n$ matrices, to combine
it with Remark~2.8 following the proposition, and to substitute the upper bound
on the expectation with a probability bound using Markov's inequality.

\begin{theorem}\label{Kostya-Liza}
 Let $\Gamma =\{\xi_{ij}\}_{1 \le i \le N, 1 \le j \le n}$ be an $N\times n$ random matrix
 on a probability space $\Omega$. Assume that entries of $\Gamma$  are independent centered
 random variables with unit variances and that in each row the entries are identically
 distributed. Let $\delta  \in (0, 1]$.
 Then there exists a random
 matrix $D_\Gamma$ on $\Omega$ taking values in $\mathcal{D}_n$ such that

 \smallskip

 \noindent
(i) for every $\omega\in \Omega$,
 $D_\Gamma(\omega)$ depends only on  the realization $\{|\xi_{ij}(\omega)|\}_{1 \le i \le N, 1 \le j \le n}$,

 \smallskip

 \noindent
(ii) for every $\omega\in \Omega$
 one has
$$
   \|R_i(\Gamma(\omega) D_\Gamma(\omega)) \|_2 \leq C\sqrt{n/\delta},
$$
(iii)  $$
  \pp\left( \mbox{\rm det } D_\Gamma \leq e^{-4 \delta N} \r)
  \leq  e^{-\delta N},
$$
where $C$ is an absolute positive constant.
\end{theorem}

As in \cite{RT}, Theorem \ref{Kostya-Liza} has important consequences.
It allows us to construct, with high probability, a diagonal matrix $D$
such that the volume of $D B_\infty^n$ remains big enough and such that,
according to Lemma~\ref{newlemma}, we have a good control of the operator
norm of $\Gamma D$ from $\ell_\infty^n$ to $\Xnk$.
Comparing to \cite{RT}, Lemma~\ref{newlemma} simplifies significantly the
proof  and allows to extend  Theorem~3.1 from \cite{RT} to the case
of rectangular matrices and to  approximations with respect to $\no_{k,2}$ norms.

\begin{theorem}\label{KL-A}
Let $1\leq n\leq N$ be integers, $\delta \in (0, 1]$. Let $k \in [N]$ such that $k \ln \frac{eN}{k} \ge n$.
Let $\Gamma$
be an $N\times n$ random matrix
satisfying the hypothesis \eqref{H}.
Then
\begin{align*}
  \pp & \left( \exists D\in \mathcal{D}_n \, | \, \mbox{ \rm det } D \geq e^{- \delta N}\, \,
  \mbox{ and } \, \,  \|\Gamma D : \ell_\infty^n \to \Xnk \| \leq C \sqrt{ \frac{kn}{\delta} \ln\left(\frac{eN}{k}\right)}  \r)
  \\
 & \geq 1- e^{-\delta N/4} -  e^{ - k \ln\left(eN/k\r)},
\end{align*}
where $C$ is a positive absolute  constant.
\end{theorem}
\proofs
Let $D_\Gamma$ be the matrix given by Theorem~\ref{Kostya-Liza}.
By property (iii) of $D_\Gamma$ it is enough to prove that
\[
 \pp\left(  \|\Gamma D : \ell_\infty^n \to \Xnk \| \leq C \sqrt{ \frac{kn}{\delta} \ln\left(\frac{eN}{k}\right)}  \r)
 \geq
 1-  e^{ - k \ln\left(eN/k\r)}.
\]

Consider two probability spaces -- the original one $(\Omega, \pp _\omega)$, where the matrix
$\Gamma$ is defined, and the auxiliary  space $(E, \pp_\eps)$, where
$E:=\{-1,1\}^{N \times n}$ and $\pp_\eps$ is the uniform probability on $E$. Given a matrix
$A =\{a_{ij}\}_{1 \le i \le N, 1 \le j \le n}$ and $\eps\in E$, denote
$A_\eps= \{\eps_{ij} a_{ij}\}_{1 \le i \le N, 1 \le j \le n}$.
Since entries of $\Gamma$ are symmetric, for every fixed $\eps\in E$ the matrix $\Gamma_\eps$ has the
same distribution on $\Omega$ as $\Gamma$.
By property~(i) of $D_\Gamma$, we have $D_\Gamma = D_{\Gamma_\eps}$ for every fixed $\eps\in E$. Therefore,
since $D_\Gamma$ is diagonal, we have for every $\eps\in E$
$$
 \left(\Gamma D_\Gamma \right)_\eps = \Gamma_\eps  D_\Gamma = \Gamma_\eps  D_{\Gamma_\eps}.
$$
Then, by  property~(ii)  of $D_\Gamma$ from Theorem~\ref{Kostya-Liza}, there exists an absolute
positive constant $C_1$ such that for every $i\leq N$ and every
$(\omega, \eps)\in \Omega\times E$,
$$
  \|R_i\left(\left(\Gamma (\omega) D_\Gamma (\omega) \right)_\eps\right) \|_2 \leq C_1\sqrt{n/\delta}.
$$
Fixing $\omega \in \Omega$ and applying  Lemma~\ref{newlemma} to the matrix $B= \Gamma(\omega) D_\Gamma(\omega)$,
we obtain that  for every fixed $\omega\in \Omega$ one has
$$
 \pp_\eps \left( \| \Gamma_\eps (\omega) D_\Gamma (\omega) : \ell_\infty^n \to \Xnk \|  > 6 C_1 \,
 \sqrt{ \frac{kn}{\delta} \ln\left(\frac{eN}{k}\right)}
\right)
\leq  e^{\displaystyle - k \ln\left(eN/k\right)}.
$$
Using that $\Gamma_\eps$ has the same distribution  as $\Gamma$ and the Fubini theorem, we obtain
\begin{align*}
  \pp_\omega & \left( \| \Gamma  D_\Gamma : \ell_\infty^n \to \Xnk  \| > 6 C_1 \,  \sqrt{ \frac{kn}{\delta} \ln\left(\frac{eN}{k}\right)}   \right)
=
\\
& =
 \pp_\eps \pp_\omega \left( \| \Gamma_\eps  (\omega) D_\Gamma (\omega) : \ell_\infty^n \to \Xnk \| > 6 C_1 \, \sqrt{ \frac{kn}{\delta} \ln\left(\frac{eN}{k}\right)}   \right)
 \\ &\leq   e^{\displaystyle - k \ln\left(eN/k\right)}.
\end{align*}
\kkk

As in Lemma~3.11 from \cite{RT}, we need to estimate the cardinality
of the set of diagonal matrices in $\mathcal{D}_n$ with not so small
determinant.
\begin{lemma} \label{KL-cardofD}
Let $n, N\geq 1$ be integers,  $\delta \in (0,1]$ and
$$
 Q:=\{D\in \mathcal{D}_n \, |\, \mbox{ \rm det } D \geq \exp(-\delta N) \}.
$$
Then  $|Q| \leq F(\delta, n, N)$, where $F(\delta, n, N)$ is defined by formula (\ref{eq:F}).
\end{lemma}

\medskip

\proofs
Note that if $D\in \mathcal{D}_n$ and $d_1$, ..., $d_n$ its diagonal elements then for every
$k\geq 0$ the set
$$
  Q_D(k) = \left\{i\leq n \,\,  |\, \, d_i = 2^{-2^{k}} \r\}
$$
has cardinality at most $m_k:=\min\{n, \lfloor {2^{-k}} 2\delta N\rfloor \}$.
Thus there are at most
$$
 \sum _{\ell =0}^{m_k} {\binom{n}{\ell}} \leq \left(\frac{en }{m_k}\r)^{m_k}
$$
choices of $\sigma _k\subset [n]$, where matrices from $\mathcal{D}_n$ may have such coordinates.
Note also that the trivial bound for the number of subsets is $2^n$.
 Denote $a:=4\delta N/n$. Note that $m_k\leq n/2$ if and only if $2^k\geq a$.

 \smallskip

\noindent
{\bf Case 1. $a\geq 2$.}
Set $ m:=\lfloor \log _2 a\rfloor \geq 1$. By above we have
\begin{align*}
  |Q|&\leq \prod_{k< m} 2^n \, \prod_{k\geq  m} \left(\frac{en }{m_k}\r)^{m_k}\leq
  2^{n m } \, \prod_{k\geq  m} \left(\frac{en }{2\delta N}\r)^{2\delta N/2^k}\,
  \prod_{k\geq  m}  2^{2k\delta N/2^k}
  \\  &\leq
  a^n \, \left(\frac{2 e }{a}\r)^{4\delta N/a}\, 2^{2\delta N (2m+1)/2^m} \leq
  (2e)^n \, a^{4\delta N/2^m}\, 2^{4\delta N/2^m} \leq (8a)^n.
\end{align*}

 \smallskip

\noindent
{\bf Case 2. $a\leq 2$.} Similarly we have
\begin{align*}
  |Q|&\leq  \prod_{k\geq  0} \left(\frac{en }{m_k}\r)^{m_k}\leq
  \prod_{k\geq  0} \left(\frac{en }{2\delta N}\r)^{2\delta N/2^k}\,
  \prod_{k\geq  0}  2^{2k\delta N/2^k}
 \leq
  \left(\frac{en }{2\delta N}\r)^{4\delta N}\, 2^{3\delta N } ,
\end{align*}
which implies the desired result.
\kkk

\subsection{Proof of Theorem \ref{Kostya}}
Let $Q$ be as in Lemma~\ref{KL-cardofD}. Note that every $D\in Q$ is diagonal with reciprocal of integers on the
diagonal. Therefore,  there exists a set ${\cal N}_D \subset T$ of cardinality
$$
|{\cal N}_D| \le N(T, \eps D B_\infty^n) \le N(T, \eps B_\infty^n)
N(B_\infty^n, D B_\infty^n)\leq M \det D^{-1} \leq M e^{\delta N}
$$
which satisfies that for any $x \in T$ there exists $y \in {\cal N}_D$ such that $x-y \in \eps D B_\infty^n$.
Let
\[
{\cal P} = \left\{ y + \eps D B_\infty^n\, | \,  D \in Q , y\in {\cal N}_D\r\}.
\]
Then, by Lemma~\ref{KL-cardofD},  $|{\cal P}| \le M e^{\delta N} F(\delta, n, N)$ and for any $x \in T$ and for any $D \in Q$ there exists $P = y_{x,D} + \eps D B_\infty^n\in \cal P$ such that $x \in  P$.

\smallskip

Theorem~\ref{KL-A} implies that with probability at least $1-e^{ - k \ln\left(eN/k\r)} -e^{-\delta N/4}$ there exists
$D\in Q$ such that
$$\Gamma(\eps D B_\infty^n)\subset C \, \eps \, \sqrt{ \frac{kn}{\delta} \ln\left(\frac{eN}{k}\right)} \, \Bnk.$$
Therefore, for such $D$,
$$
  \Gamma (x-y_{x,D})  \in \Gamma ( \eps D B_\infty^n) \subset  C \, \eps \,
   \sqrt{ \frac{kn}{\delta} \ln\left(\frac{eN}{k}\right)} \, \Bnk,
$$
hence,
$$
  \Gamma (P) \subset \Gamma x + \Gamma (y_{x,D}-x) +\Gamma ( \eps D B_\infty^n)
  \subset
  \Gamma x + 2 C\, \eps \, \sqrt{ \frac{kn}{\delta} \ln\left(\frac{eN}{k}\right)} \, \Bnk.
$$
This proves the existence of a ``good" collection $\mathcal{P}$.

Finally, let
$\mathcal{P}'$ be the set of all $P\in \mathcal{P}$  such that $P\cap T\ne \emptyset$.
For every $P\in \mathcal{P}'$  choose an arbitrary $z_P\in P\cap T$ and let $\mathcal{N}=\{z_P\}_{P\in \mathcal{P}'}$.
 By above, for every $x\in T$ there exists $D\in Q$ and $P=y_{x,D} + \eps D B_\infty^n \in \mathcal{P}$
 such that $x\in P$, in particular $P \in \mathcal{P}'$, and
$$
 \Gamma (P) \subset   \Gamma x + 2 C \, \eps \, \sqrt{ \frac{kn}{\delta} \ln\left(\frac{eN}{k}\right)} \, \Bnk.
$$
Thus,
$\Gamma z_P \in \Gamma x  +2C \, \eps \, \sqrt{ \frac{kn}{\delta} \ln\left(\frac{eN}{k}\right)} \, \Bnk.$ This implies the desired  result.
\kkk

\begin{remark}
\label{rem:appli}{\rm
We apply Theorem \ref{Kostya} for $T\subset t B_2^n$, $t\geq 1$, $\eps \leq 1/\sqrt n$, and
$\delta \geq n/(2N)$ so that, $F(\delta, n, N) = \left(32\delta N/n\r)^n$.
 Then Theorem~\ref{Kostya} combined with
Lemma~\ref{covcubesch} implies that
there exists $\mathcal{N}\subset T$ with cardinality at most
$$
  \left(\frac{224\delta t N}{\eps n^{3/2}}\r)^n\,  e^{\delta N}
$$
such that with probability at least $1- e^{- k \ln\left({eN}/{k}\right)}-e^{-\delta N/4}$ one has
$$
 \forall x\in T \,\, \exists y\in \mathcal{N} \, \, \, \, \mbox{ such that } \, \, \, \,
  \Gamma (x-  y) \in  C \, \eps \, \sqrt{ \frac{kn}{\delta} \ln\left(\frac{eN}{k}\right)} \, \Bnk.
$$
}
\end{remark}

\section{Geometry of Random Polytopes}
\label{geom}
In this section, we study some classical geometric parameters associated to  random
polytopes of the form $K_N:= \Gamma ^* B_1^N$, where $\Gamma = \{\xi_{ij}\}_{1 \le i \le N, 1 \le j \le n}$
is an $N\times n$
random matrix. In other words, $K_N$ is the
absolute convex hull of the rows of $\Gamma$.
We provide estimates  on the asymptotic behavior of the volume and
the mean widths of $K_N$ and its polar.
In this section, the random operator $\Gamma$ satisfies the hypothesis \eqref{H}:
the random variables $\xi_{ij}$ are independent symmetric with unit variances such
that in each row of $\Gamma$ the entries are identically distributed. Moreover,
we assume that the random variables  $\xi_{ij}$  satisfy a uniform small ball
probability condition which means that we can fix $u, v \in (0,1)$ such that
\[
\forall i,j \quad\sup\limits_{\lambda\in \RR}\pp\bigl\{ |\xi_{ij}-\lambda|\leq u\bigr\} \leq v.
\]

\subsection{Inclusion Theorem}
We start  by showing that for an
$N\times n$ random matrix $\Gamma $ satisfying conditions described above,
the body $K_N = \Gamma ^* B_1^N$
contains a large ``regular'' body  with high probability.

\begin{theorem} \label{inclusion}
Let $\beta \in (0, 1)$.
There are two positive constants $\KK=\KK(u, v, \beta)$ and $C(u, v, \beta)$
which depend only on $u, v, \beta$ and an absolute constant $c>0$, such that the following holds.
For every positive integers $n, N$ satisfying $N\geq Mn$ one has
$$
\pp \left( K_N \supset C(u, v, \beta)
\left( B_{\infty}^n \cap \sqrt{\ln (N/n)}\,  B_2^n
\r) \r) \geq 1 -  4 \exp \left( -c n^{\beta} N^{1-\beta} \r).
$$
\end{theorem}

\medskip

\begin{remark}
{\rm
 It is known that for a  Gaussian random matrix
one has
$$
\pp \left( K_N \supset C\sqrt {\beta \ln (N/n)} \ B_2^n  \r)
  \geq 1- 3 \exp
\left( -c n^{\beta} N^{1-\beta} \r),
$$
where $C$, $c$ are  absolute positive constants (see e.g. \cite{Gl}).
Moreover, the probability estimate cannot be improved.
Indeed, for a Gaussian random matrix and $\beta \in (0, c'')$ one has
$$
\pp \left( K_N \supset   C'\sqrt {\beta \ln (N/n)} B_2^n  \r)
\leq 1-\exp   \left( -c' n^{\beta} N^{1-\beta} \r),
$$
where $C', c' >0$ and  $ 0 < c'' \le 1$ are  absolute constants.
}
\end{remark}
Since $B_{\infty}^n \subset \sqrt{n} B_2^n$, Theorem~\ref{inclusion} has
the following consequence.
\begin{sled} \label{maincor}
 Under the assumptions and notations of Theorem~\ref{inclusion}, for  $\KK n<N\leq e^n$  one has
$$
\pp \left( K_N\supset C(u, v, \beta) \sqrt{\frac{ \ln
( N/n)}{n}} \, B_{\infty}^n \r)   \geq 1 -  4 \exp \left( -c n^{\beta} N^{1-\beta} \r).
$$
\end{sled}

In fact, our proof of Theorem \ref{inclusion} gives that if
$$
   N \geq n \max\left\{ \exp (4 C_{v}/ \beta) ,
\left(\frac{C \ln (e/(1-\beta)}{c_{uv}\, (1-\beta)}\r)^{1/(1-\beta)} \r\},
$$
where $C > 1$ is an absolute positive constant,
$c_{uv}=c u v \sqrt{1- v }$ is the constant from Lemma~\ref{levy} below,
 and
$C_{v} =  5 \ln (2/(1-v))$,
 then
\begin{equation}
\label{eq:inclusion-precise}
\pp \left( K_N \supset \frac{c_{uv}}{2 \sqrt 2}\,
\left( B_{\infty}^n \cap R B_2^n
\r) \r) \geq 1 -  4 \exp \left( -\frac{n^{\beta} N^{1-\beta}}{40} \r)
\end{equation}
with $R=\sqrt{\beta \ln (N/n) /C_{v}}$.
Note that
$K_N = \mbox{abs conv} \{x_j \} _{j\leq N}$,
where $x_j = \Gamma^* e_j$ are the columns of $\Gamma^*$.
Hence for every $z\in \R$,
\[
h_{K_N} (z) = \sup _{j\leq N} |\la z, x_j\ra| =\| \Gamma z\|_{\infty}.
\]
Let
$L= c_{uv} (B_{\infty}^n \cap R B_2^n)$.
To prove \eqref{eq:inclusion-precise}, we show that
\begin{equation}
\label{eq:mainprob}
\pp \left( \exists z \in \partial L^0\,\,\, |\,\,\, \| \Gamma z \|_\infty <  \frac{1}{4} \r)
\le
4 \exp \left( -\frac{n^{\beta} N^{1-\beta}}{40} \r).
\end{equation}
The proof of this statement will be divided into two steps.
First, we will show an individual estimate for a fixed
$z \in \partial L^0$. Then we
use the
net introduced in Theorem~\ref{Kostya} to get a global
estimate for any point of this net, using that  this net is a subset
of $\partial L^0$. A crucial point is that this net is a good
covering of $\Gamma (\partial L^0)$ in $\nk\cdot$-metric.

\subsubsection{Basic facts about small ball probabilities.}
Recall that for a (real) random variable $\xi$ its {\it L\'evy concentration function}
$\cf(\xi, \cdot)$ is defined on $(0, \infty)$ as
$$
 \cf(\xi, t):=\sup\limits_{\lambda\in \RR}\pp\bigl\{ |\xi-\lambda|\leq t\bigr\}.
$$
For any centered random variable with unit variance, there exist $u, v \in (0,1)$ such that
\begin{equation}\label{matr-cond-conc}
  \cf(\xi, u) \leq v.
\end{equation}
The following lemma is a consequence of  Rogozin's theorem  \cite{ROG}
that was used for example in \cite{RT} (see Lemma~4.7 there).
\begin{lemma}\label{levy}
Let $\xi_1, ..., \xi_m$ be independent random variables   satisfying
(\ref{matr-cond-conc}) with the same $u, v \in (0,1)$. Then for every $x\in S^{m-1}$ one has
$$
  \cf\big(\sum _{i=1}^m x_i \xi_i, c_{uv} \big)\leq v,
$$
where $c_{uv}=c u v \sqrt{1- v }$ and $c\in (0,1]$ is an absolute constant.
\end{lemma}
\begin{remark}{\rm
If we have a bounded  moment of order larger than 2, then we could use
a consequence of the Paley-Zygmund inequality, which also provides a lower bound on the
small ball probability of a random sum. The following statement was proved in \cite[Lemma~3.1]{LR}
following the lines of \cite[Lemma~3.6]{LPRT} with appropriate modifications to deal with centered
random variables (rather than symmetric):}

Let $2<r \leq 3$ and $\mu \geq 1$.
Suppose $\xi_1, \ldots, \xi_m$ are independent centered random variables
such that $\E |\xi_{i}|^2\geq 1$ and $\E |\xi_{i}|^r \leq \mu^r$ for every $i \leq m$.
Let $x = (x_i) \in \ell_{2}$ be such that $\|x\|_2 = 1$.
Then for every $\lambda \geq 0$
\begin{equation}\label{omar}
\P \left( \Big| \sum_{i=1}^{m} \xi_i x_i \Big| > \lambda \r) \geq
\left( \frac{ 1 - \lambda^2 }{ 8\mu^2 } \right)^{r/(r-2)}.
\end{equation}
\end{remark}
%
%

\noindent
{\bf Proof of Lemma~\ref{levy}. }
Fix $x\in S^{m-1}$. We clearly have
$\cf(x_i \xi_i, |x_i| u \big)\leq v$ for every $x_i\ne 0$.
Applying Theorem~1 of \cite{ROG} to random variables $x_i \xi _i$, $i\leq m$,
we observe there exists and absolute constant $C\geq 1$ such
that for every $w\geq u\|x\|_\infty /2$,
$$
  \cf\big(\sum _{i=1}^m x_i \xi_i, w \big)\leq
 \frac{ C w}{\sqrt{ \sum _{i=1}^m  |x_i|^2 u^2 \big(1-\cf(x_i \xi_i,|x_i| u ) \big)} }
 \leq
 \frac{ C w}{ u \sqrt{1- v } }.
$$
Take $w=  u v \sqrt{1- v }/C$. If $\|x\|_\infty \leq 2 v \sqrt{1- v }/C$
then $w\geq u\|x\|_\infty /2$. Therefore for such $x$ we have
$$
  \cf\bigg(\sum _{i=1}^m x_i \xi_i, w \bigg)\leq v.
$$
If there exists $\ell \leq m$ such that $|x_\ell| > 2 v \sqrt{1- v }/C$, then we have
$$
  \cf\bigg(\sum _{i=1}^m x_i \xi_i, w \bigg)\leq \cf\big(x_\ell \xi_\ell, w \big) =
  \cf\big( \xi_\ell, w/ |x_\ell| \big) \leq \cf\big( \xi_\ell, u \big) \leq  v,
$$
which completes the proof.
\kkk

\subsubsection{The individual small ball estimate.}
To prove Theorem~\ref{inclusion} we need to extend a result by
Montgomery-Smith \cite{MS}, which originally was proved for Rademacher
random variables. Note that this lemma does not require any conditions on
the moments of random variables.
%
%
%
\begin{lemma} \label{monsm}
Let $\xi _i$, $i\leq n$, be independent symmetric random
variables satisfying  condition (\ref{matr-cond-conc}). Let
$\alpha \ge 1$ and
$L= c_{uv} (B_{\infty}^n \cap \alpha B_2^n)$,
where  $c_{uv}$ is the constant from Lemma~\ref{levy}.
Then for every non-zero $z\in \R$ one has
$$
 \pp \left( \sum_{i=1}^n \xi _i z_i   > h _L (z) \r)
 > ((1-v)/2)^{5\alpha^2} .
$$
\end{lemma}

We postpone the proof of this lemma to the end of this section. Note that if our variables satisfy
$1\leq \E \xi_i ^2 \leq \E |\xi_i |^r \leq \mu ^r$ for some $r > 2$ then using (\ref{omar}) and
repeating the proof of Lemma~4.3 from \cite{LPRT} we could consider
$L= (1-\delta) (B_{\infty}^n \cap \alpha B_2^n)$ and estimate the corresponding
probability from below by $\exp \left( -C_{\mu, \delta, r} \alpha ^2   \r)$, where
$C_{\mu, \delta, r}$ depends only on  $\mu, \delta, r$.

\medskip

Lemma~\ref{monsm} has the following consequence.
\begin{lemma}
\label{cmonsm}
Under assumptions of Lemma~\ref{monsm}
for every $z\in \R$ and  every
$\sigma \subset [N]$ one has
$$
  \pp \left( \| P_{\sigma} \Gamma z\| _{\infty}
   < h_L (z) \r) <
   \exp\left( -|\sigma| \exp (- C_{v} \alpha ^2)\r),
$$
where $P_{\sigma} : \ \RR^N \to \RR ^{\sigma}$ is the coordinate projection
and $C_v=5\ln (2/(1-v)))$.
\end{lemma}

\proofs
Applying  Lemma~\ref{monsm} to the $|\sigma|\times n$ random matrix
$P_{\sigma} \Gamma = (\xi _{i j}) _{i\in \sigma, j\leq n}$
 we have for every
$z=\{z_j\}_{j=1}^n \in \R$ and every $i\in \sigma $
$$
\pp \left( \sum _{j=1}^n z_j \xi _{i j}
< h_L (z) \r) \leq 1-  \exp( -C_v \alpha^2 )\leq
\exp\left( - \exp( -C_v \alpha^2) \r).
$$
Thus
\begin{eqnarray*}
\lefteqn{\pp \left( \| P_{\sigma} \Gamma z\| _{\infty}
< h_L (z) \r) =
\pp \left( \sup _{i \in \sigma}
\left| \sum _{j=1}^n z_j \xi _{i j} \r|
< h_L (z) \r)}\\
& = & \prod _{i \in \sigma} \pp
\left( \left| \sum _{j=1}^n z_j \xi _{i j} \r|
< h_L (z) \r) <
\exp\left(-|\sigma| \exp (-C_{v} \alpha ^2)\r).
\end{eqnarray*}
\kkk
We can now state the main individual small ball estimate. \begin{lemma}
\label{lem:indiv}
Let $\beta \in (0,1)$ and define $m = 8 \lceil (N/n) ^{\beta} \rceil$
(if the latter number is larger than or equal to $N/4$ we take $m=N$) and $k=\lfloor N/m\rfloor$.
Let $L= c_{uv} (B_{\infty}^n \cap R B_2^n)$, where $R=\sqrt{\beta \ln (N/n) /C_{v}}$.
Then for any $z \in \partial L^o$ one has
\[
\pp \left(\frac{1}{\sqrt k } \nk{\Gamma z} <  \frac{1}{2} \r)
\le
\exp(- 0.3  \ n^\beta \, N^{1-\beta}).
\]
\end{lemma}
\proofs
Below we assume $m<N/4$ (then $k\geq 4$, hence
$km > 4N/5$); the proof in
the case $m=N$, $k=1$ repeats the same lines with simpler
calculations. Let $\sigma _1$, \ldots, $\sigma _k$ be a
partition of $[N]$ such that
$m\leq |\sigma _i|$ for every $i\leq k$.
Then, for any $a \in \RR^N$
\[
\frac{1}{\sqrt k} \nk{a}
\ge
\frac{1}{\sqrt k} \left( \sum_{i=1}^k \| P_i a \| _{\infty}^2 \r)^{1/2}
\ge \frac{1}{k} \sum _{i=1}^k \| P_i z \| _{\infty},
\]
where $P_i = P_{\sigma _i} :  \RR^N \to \RR ^{\sigma _i}$
is the coordinate projection.
Define $|||\cdot |||$ on $\RR^N$ by
$$
||| z ||| = \frac{1}{k} \sum _{i=1}^k \| P_i z \| _{\infty}
$$
for every $z\in \RR^N$.
Note that if for some $z\in \R$ we have
$||| \Gamma z ||| < h_L (z) /2$ then there exists
$I\subset [k]$ of cardinality
at least $k/2$ such that for every $i\in I$ one has
$\| P_i \Gamma z \| _{\infty} < h_L (z) $.
Applying Lemma~\ref{cmonsm} with  $\alpha = R$
(note that $\alpha \ge 2$, by the condition on $n$ and $N$),
we obtain
for every $z=\{z_i\}_{i=1}^n \in \R$,
\begin{align*}
{\pp \left( |||\Gamma z ||| < h_L (z) /2 \r)}
& \le    \sum _{|I|=[(k+1)/2]}
\pp \left( \| P_i \Gamma z \| _{\infty} < h_L (z)
\, \mbox{ for every } \, i \in I \r)
\\
& \le  \sum _{|I|=[(k+1)/2]}
\prod _{i\in I}\pp
\left( \| P_i \Gamma z \| _{\infty} < h_L (z) \r)\\
& \le  \sum _{|I|=[(k+1)/2]} \prod _{i\in I}
\exp\left(-|\sigma _i| \exp (-C_{v} \alpha ^2)\r)\\
& \le  \binom{k}{[k/2]}  \exp\left(- (km/2)
\exp (-C_{v} \alpha ^2 )\r) \\
& \le   \exp\left(
k \ln 2 - (km/2) \exp (-C_{v} \alpha ^2)\r),
\end{align*}
where $C_{v} =  5 \ln (2/(1-v))$.
By our choice of $k$ and $m$ we have $km >4N/5$, therefore
$(km/2) \exp (-C_{v} \alpha ^2)\geq 2 N^{1-\beta}  n^{\beta} /5$.
We also have $k\leq N^{1-\beta} n^{\beta} /8$. Thus
$$
\pp \left( |||\Gamma z ||| < h_L (z) /2 \r) \le
\exp\left( -0.3\  N^{1-\beta} n^{\beta} \r).
$$
This completes  the proof.
\kkk

 Finally we prove  Lemma~\ref{monsm}.
 For a positive integer $m$, define $||| \cdot |||_m$  on $\R$ by
\begin{equation*} 
 ||| z |||_m = \sup \sum _{i=1}^m \left(\sum _{k\in B_i}
 |z_k |^2  \r)^{1/2},
\end{equation*}
 where the supremum is taken over all partitions
$B_1, \ldots, B_m$ of $[n]$. We will need the following lemma,
which was essentially proved in \cite{MS} (see Lemma~2 there).
\begin{lemma}\label{ms}
Let $\alpha \geq 1$ and $m\geq 1+4\alpha^2$ be an integer. For all $x\in \R$ one has
$$
  h_{ B_\infty^n\cap \alpha B_2^n}(x) \leq |||x|||_m.
$$
\end{lemma}

\smallskip

\proofs
Fix $x\in \R$ and choose
$y\in B_\infty^n\cap \alpha B_2^n$ so that $h(x) = \sum _i x_i y_i$. For every $k$ with
$y_k^2 \geq 1/2$ choose $B_{1, k}=\{k \}$. Since $|y| \leq \alpha $
there are at most  $2\alpha^2 $ such sets. Denote $B:=\cup_k  B_{1,k}$.
Now let $z_i$ denote $y_i$ if $|y_i|\leq 1/\sqrt 2$ and
$z_i=0$ otherwise. Let $n_0=0$ and define $n_0< n_1< n_2<...$ by
$$
   n_{k+1}= 1 +\sup \big\{ \ell \in  [n_k+1, n-1]\,\,  | \, \sum _{i=n_k+1}^\ell z_i^2\leq  1/2\big\}
$$
(if $n_k=n$ we stop the procedure). Denote $B_{2, k}:=[n_{k-1}+1, n_k]\setminus B$.
Since $|y| \leq \alpha $ we have at most $2\alpha^2 +1$ such sets. Moreover, we have
 $$\sum _{i\in B_{2, k}} z_i^2 = \sum _{i\in B_{2, k}} y_i^2 \leq 1.$$ Since $y\in B_\infty^n$ and
 $m\geq 4\alpha^2 +1$,  we obtain
$$
    h(x) = \sum _{i=1}^n x_i y_i \leq \sum _{j=1}^2 \sum_k
    \bigg(\sum _{i\in B_{j, k}}  x_i^2\bigg)^{1/2}\bigg(\sum _{i\in B_{j, k}}  y_i^2\bigg)^{1/2}
    \leq  \sum _{j\leq 2, k}
    \bigg(\sum _{i\in B_{j, k}}  x_i^2\bigg)^{1/2} \leq |||x|||_m.
$$
\kkk

\noindent
{\bf Proof of Lemma~\ref{monsm}: }
We folow the lines of  Montgomery-Smith's proof.
Let $m=\lceil 1+4\alpha^2\rceil$.
Given $z \in \R$, let $m' \le m$ and
$B_1, \ldots, B_{m'}$ be a partition of  $[n]$
such that
$$
\forall i\leq m' \, \, \, \sum _{k\in B_i} |z_k |^2 \ne 0 \quad \quad
\mbox{ and } \quad \quad  ||| z |||_m =  \sum _{i=1}^{m'}
\left(\sum _{k\in B_i} |z_k |^2 \r)^{1/2}.
$$
Then, using Lemma~\ref{ms}, we have
\begin{eqnarray*}
\lefteqn{p:= \pp \left( \sum _{i=1}^n \xi _i z_i > h_L (z) \r) \geq
\pp \left( \sum _{i=1}^n \xi _i z_i > c_{uv} \, ||| z |||_m  \r) }\\
& = & \pp \left( \sum _{i=1}^{m'} \sum _{k\in B_i} \xi _k z_k >
c_{uv} \, \sum _{i=1}^{m'} \left(\sum _{k\in B_i}
| z_k |^2 \r)^{1/2}\r)\\
& \ge & \pp \left(\bigcap_{i \le m'}\left(\sum _{k\in B_i}
\xi_k z_k \ge
c_{uv} \, (\sum _{k\in B_i} |z_k|^2)^{1/2}\right) \right).
\end{eqnarray*}
Since $\xi _i$'s are independent we obtain
$$
p \geq \prod _{i=1}^{m'} \pp \left( \sum _{k\in B_i} \xi _k z_k
> c_{uv} \,  \left(\sum _{k\in B_i} |z_k |^2
\r)^{1/2}\r).
$$
For $i \le m'$ set
$$
f_i =  \left( \sum _{k\in B_i} \xi _k z_k \r) \cdot
\left(\sum _{k\in B_i} |z_k |^2 \r)^{-1/2} .
$$
Using that $\xi _i$'s are symmetric and applying Lemma~\ref{levy}
we get
\begin{eqnarray*}
\pp \left( f_i >  c_{uv}  \r) & = &
\frac{1}{2}\ \pp \left( |f_i| >  c_{uv}  \r)
 \geq  \frac{1-v}{2}.
\end{eqnarray*}
Therefore,
$$
  p \geq ((1-v)/2)^{m'}\geq ((1-v)/2)^{m} \geq ((1-v)/2)^{5\alpha^2},
$$
which implies the desired result.
\kkk
\subsubsection{The global small ball estimate.}
In this section, we prove Theorem~\ref{inclusion}.
As we mentioned after its statement, our goal is to
prove \eqref{eq:mainprob} for $N \ge M n$,
where $M$ depends only on $\beta, u$ and $v$.

Let $\beta \in (0,1)$ and, as in Lemma~\ref{lem:indiv}, define $m = 8 \lceil (N/n) ^{\beta} \rceil$  and $k=\lfloor N/m\rfloor$ so that $N^{1-\beta} n^{\beta} / 10 \le k \le N^{1-\beta} n^{\beta} / 8$.  By the choice of $M$, we obviously have $k \ln (eN/k) \ge n$.
Let $T = \partial L^o$ and set
$$
  \delta = 0.1 (n/N)^{\beta}  \quad  \mbox{ and }  \quad
  \eps= \frac{1}{c_{uv}\, \sqrt{n} \, \exp( (N/n)^{1-\beta}/20)} .
$$
 Since
$$
 T\subset L^0=c_{uv}^{-1}\big(\conv B_1^n\cup (B_2^n/R)\big) \subset c_{uv}^{-1} B_2^n,
$$
 we use Theorem~\ref{Kostya} (see Remark~\ref{rem:appli})
  to construct  a set $\mathcal{N} \subset T$ of cardinality at most
$$
  \left(\frac{224\delta N}{\eps c_{uv} n^{3/2}}\r)^n\,  e^{\delta N}
$$
such that with probability at least $1- e^{- k \ln\left({eN}/{k}\right)}-e^{-\delta N/4}$ one has
\begin{equation}\label{k-l}
  \forall  x\in T \, \, \exists z\in \mathcal{N} \quad \mbox{ such that }
  \quad \nk{\Gamma(x-z)} \leq    C_1 \, \eps \, \sqrt{ \frac{kn}{\delta} \ln\left(\frac{eN}{k}\right)},
\end{equation}
where $C_1>0$ is an absolute constant.
Since
$$
\exp\left(  n \ln (224 \delta N/n)  +
            n \ln (1/(\eps c_{uv} n^{1/2}) ) + \delta N  -0.3\  N^{1-\beta} n^{\beta} \r)
 \le
 \exp\left(  -0.1\  N^{1-\beta} n^{\beta} \r),
$$
provided that $(N/n)^{1-\beta}$ is large enough, and ${\cal N } \subset T$, we deduce from Lemma~\ref{lem:indiv} that
\begin{align*}
\pp    \left(  \exists z\in \mathcal{N} \, : \, \frac{1}{\sqrt k} \nk{\Gamma z} < 1/ 2
\r)
\le  \sum _{z\in\mathcal{N} } \pp \left( \frac{1}{\sqrt k} \nk{\Gamma z} < 1/ 2 \r)
 \le
 \exp\left(  -0.1\  N^{1-\beta} n^{\beta} \r).
\end{align*}
Let $\overline{\Omega}$ be the subset of $\Omega$, where (\ref{k-l}) holds.
Then,
on $\overline{\Omega}$, for every $x\in T$
there exists $z\in \mathcal{N}$ such that
\begin{align*}
\frac{1}{\sqrt k} \nk{\Gamma z} & \leq \frac{1}{\sqrt k} \nk{\Gamma x} + \frac{1}{\sqrt k} \nk{\Gamma (z - x)}
\leq \frac{1}{\sqrt k} \nk{\Gamma x} +  C_1 \eps \, \sqrt{ \frac{n}{\delta} \ln\left(\frac{eN}{k}\right)}.
\\
& \le \frac{1}{\sqrt k} \nk{\Gamma x} + \frac{C_2 \sqrt{\left(\frac{N}{n}\r)^{\beta} \ln \left(10 e \left(\frac{N}{n}\right)^\beta\r)}}{c_{uv}\, \exp((N/n)^{1-\beta}/20) } ,
\end{align*}
where $C_2$ is an absolute positive constant. Since $N \ge \KK n$
(for large enough $\KK$ depending only on $u, v$ and $\beta$), we observe
\[
c_{uv}^2 \exp((N/n)^{1-\beta}/10)> 16 \, C_2^2 \left(\frac{N}{n}\r)^{\beta}
\ln \left(10 e \left(\frac{N}{n}\right)^\beta\r).
\]
Therefore,
\begin{align*}
 \pp & \left( \left\{ \omega \in \overline{\Omega} \ | \
\exists x\in \partial L^o \, : \, \frac{1}{\sqrt k} \nk{\Gamma x} < \frac{1}{4} \r\} \r)
 \\ &\leq
   \pp \left( \left\{ \omega \in \overline{\Omega} \ | \
\exists z\in \mathcal{N} \, : \, \frac{1}{\sqrt k} \nk{\Gamma z} < \frac{1}{ 2} \r\} \r)  \leq
\exp\left(  -0.1\  N^{1-\beta} n^{\beta} \r)  .
\end{align*}
The desired result follows since $h_{K_N} (x) = \| \Gamma x \|_\infty \geq \frac{1}{\sqrt k} \nk{\Gamma x} $
for every $x\in \R$ and since
$$\pp (\overline{\Omega} )\geq 1- e^{- k \ln\left({eN}/{k}\right)}-e^{-\delta N/4} \ge 1-2 \exp(N^{1-\beta} n^{\beta}/40).
$$
\kkk

\subsection{Volumes and mean widths of $K_N$ and $K_N^0$}

In this section we apply the results of the
previous subsection to obtain asymptotically sharp
estimates for the volumes and  the
mean widths of $K_N$ and  $K_N^0$. We refer to \cite{P}
for general knowledge about these parameters.
We recall that by  Santal\'{o} inequality and Bourgain-Milman
\cite{BM} inverse Santal\'{o} inequality there exists an
absolute positive constant $c$ such that for every convex
symmetric body $K$ one has
\begin{equation} \label{sant}
c^n |B_2^n| ^2 \leq |K| |K^0| \leq |B_2^n| ^2 .
\end{equation}
Below we fix constants $\KK=\KK(u, v, \beta)$ and $C(u, v, \beta)$ from Theorem~\ref{inclusion}.

We start estimating the volumes of $K_N$ and $K_N^0$.
For convenience we separate upper and lower
estimates (some bounds require an additional condition on the matrix $\Gamma)$.
Corollary~\ref{maincor} and (\ref{sant})
 imply the following volume
estimates for $K_N$ and $K_N^0$.

\begin{theorem}
\label{volone}
Let $\KK n<N\leq e^n$, $\beta\in (0, 1)$.
There exists  absolute positive constants $C$ and $c$ such that
with probability at least  $1 -  \exp \left( -c n^{\beta} N^{1-\beta} \r)$ one has
$$
|K_N|^{1/n} \geq 2 C(u, v, \beta)
       \sqrt{\frac{ \ln (N/n)}{n}} \, \, \, \,\, \,\, \,
\mbox{ and } \, \, \, \,\, \,\, \, |K_N^0|^{1/n} \leq
\frac{C}{C(u, v, \beta)  \sqrt{  n \ln (N/n)}}.
$$
\end{theorem}

\medskip

To prove the remaining bounds on volumes of $K_N$ and $K_N^0$ we
introduce one more condition on the matrix $\Gamma$, namely we require that
\begin{equation}\label{cond-two}
 \P\left(\max_{i\leq N}  |\Gamma ^* e_i| > \lambda  \sqrt{ n} \r) \leq p_0
\end{equation}
for some $0<p_0<1$ and $\lam \geq 1$. Such condition holds for example when entries of $\Gamma$
are i.i.d. centered random variables with finite $p$-th moment for some $p>4$,
provided that $N\leq C_p n^{p/4}$ (this can be proved using Rosenthal's inequality,
see Corollary~6.4 in \cite{GLPT}).

\smallskip

The lower bound on $|K_N|$ (and the upper bound on $|K_N^0|$)
follows from (\ref{sant}) and a well known estimate on the volume of the convex hull of $k$
points (\cite{BF}, \cite{CP}, \cite{Gl}):

\smallskip

{\it Let $2n\leq k \leq e^n$ and $z_1, \ldots, z_k \in S^{n-1}$, then
\begin{equation*}
|\mbox{\rm abs conv} \{z_i\}_{i\leq k} | ^{1/n} \leq c
\sqrt{\ln (k/n)}/ n,
\end{equation*}
where $c>0$ is an absolute  constant.}

\begin{theorem} \label{voltwo} Let $\KK n<N\leq e^n$ and $\beta \in (0, 1)$.
Assume that the matrix  $\Gamma$ satisfies (\ref{cond-two}).
There exist absolute positive constants $c$ and $C$ such that
one has
$$
|K_N|^{1/n} \leq C \lambda \sqrt{\frac{\ln (N/n)}{n}} \, \, \, \,\, \,\, \,
\mbox{ and } \, \, \, \,\, \,\, \, |K_N^0|^{1/n} \ge
c/(\lambda \sqrt{n \ln (N/n)})
$$
with probability at least $1-p_0$.
\end{theorem}

An important geometric parameter associated to a convex body is the (half of) mean width
of $K^0$ defined by
$$
M_K=M(K)=\int _{S^{n-1}} \|x \| _K \ d\nu,
$$
where $\nu$ is the normalized Lebesgue measure on $S^{n-1}$.  It is well known that there exists a
constant $c_n >1$ ($c_n \to  1$ as $n\to \infty $) such that
$$
 M_K = \frac{c_n}{\sqrt{n}} \E \| \sum_{i=1}^{n} e_i g_i\|_K,
$$
for every $K \subset \R$.
The (half of) mean width of $K$, $M(K^0)$, we denote by $M^*_K = M^*(K)$.  Observe that
\[
M^*(K) =  \frac{c_n}{\sqrt{n}} \E \| \sum_{i=1}^{n} e_i g_i\|_{K^0} =  \frac{c_n}{\sqrt{n}} \E \sup_{t \in K} \sum_{i=1}^{n} t_i g_i =  \frac{c_n}{\sqrt{n}} \ell_*(K),
\]
where $\ell_*(K)  = \E \sup_{t \in K} \sum_{i=1}^{n} t_i g_i $ is
the  Gaussian complexity measure of the convex body $K$.
We recall the following inequality, which holds for every
convex body $K$ (see e.g. \cite{P})
\begin{equation} \label{ury}
M_K^* \geq \left(|K| / |B_2^n| \r) ^{1/n} \geq 1/M_K .
\end{equation}
%
Now we calculate the mean widths $M(K_N)$ and
$M(K_N^0)$.
\begin{theorem} \label{mk}
Let $\KK  n<N\leq e^n$ and $\beta \in (0, 1)$. Then
$$
 M(K_N) \leq C C^{-1}(u, v, \beta) \left( \sqrt{(\ln (2n) )/n }+1/ \sqrt{ \ln(N/n)}  \r)
$$
with probability at least
$1 -  \exp \left( -c n^{\beta} N^{1-\beta} \r),$
where $C$ and $c$ are  absolute positive constants.
Moreover, if  the matrix  $\Gamma$ satisfies (\ref{cond-two}), then
there exists an absolute positive constant $c_1$ such that
with probability at least $1-p_0$ one has
$$
 M(K_N) \geq c_1 /(\lambda  \sqrt{\ln (N/n)}).
$$
\end{theorem}

\proofs
By Theorem~\ref{inclusion} we have
\begin{eqnarray*}
M(K_N) & \leq & M\left( C(u, v, \beta) \left( B_{\infty}^n \cap
\sqrt{ \ln(N/n)} B_2^n \r) \r) \\
& \leq   & \left(1/C(u, v, \beta)\right) \left( M \left( B_{\infty}^n \r) +
M\left( \sqrt{\ln(N/n)} B_2^n \r)
\r) ,
\end{eqnarray*}
which  proves the upper bound.

By (\ref{ury}) and Theorem~\ref{voltwo}  there exists
an absolute positive constant $c_1$ such that
$$
M(K_N) \geq \left( |B_2^n|/|K_N|\r)^{1/n}
\geq c_1 /(\lambda  \sqrt{\ln (2N/n)}),
$$
with probability larger than or equal to $1-p_0$.
This proves the lower bound.
\kkk

\begin{remark}{\rm
 Note that by Theorem~\ref{mk},
for $N\le \exp(n/\ln n)$ we have
$$
M(K_N) \approx 1 / \sqrt{\ln (N/n)}.
$$
If $N\geq \exp(n/\ln(2n))$ there is a gap between
lower and upper estimates. Both estimates could be
asymptotically sharp as was shown in \cite{LPRT}.
}
\end{remark}


\begin{theorem}
\label{mstar}
There exist positive absolute constants $c$,  $c_0$, and $C$
such that the following holds.
Let $\KK n<N\leq e^n$. Then
$$
  M(K_N^0) \geq c_0 \sqrt{\ln (N/n)}
$$
with probability at least $1 -  \exp \left( -c n^{\beta} N^{1-\beta} \r)$.
Moreover, assuming that
 the matrix  $\Gamma$ satisfies (\ref{cond-two}), with probability at least $1-p_0$
 one has
$$
 M (K_N^0) \leq C \lambda  \sqrt{\ln N}.
$$
\end{theorem}

\proofs
By (\ref{ury}) we have
$$
 M(K_N^0)\geq \left(|B_2^n|/|K_N^0|\r)^{1/n} .
$$
Therefore, the lower bound follows by Theorem~\ref{volone}.

Let $G=\sum_{i=1}^n g_i e_i$. Recall that $K_N$ is the absolute
convex hull of $N$ vertices $\Gamma ^* e_i$. Thus we have
$$
 M(K_N^0) \leq \frac{c_1}{\sqrt{n}} \E  \| G \| _{K_N^0}
 = \frac{c_1}{\sqrt{n}} \E \max _{i\leq N} \la G, \Gamma ^* e_i \ra ,
$$
where $c_1$ is an absolute constant. Since with probability at least
$1-p_0$ we have $|\Gamma ^* e_i| \leq \lambda \sqrt{n}$ for every
$i\leq N$, using standard estimate for the expectation of maximum
of Gaussian random variables (see, e.g., \cite{P}), we obtain that
there is an absolute constant $c_2$ such that
$$
M(K_N^0) \leq c_2 \lam \ \sqrt{\ln N}
$$
with probability larger than or equal to $1-e^{-n}$.
\kkk

Finally we note that the bounds of Theorem~\ref{mstar} are sharp, whenever
$\ln N$ and $\ln (N/n)$ are comparable, for example if $N>n^2$. However, when
$N$ is close to $n$ we have a gap between upper and lower bounds. Below we provide
a better lower bound for $ M(K_N^0)$ in the case $N\leq n^2$, which closes this gap.
We will need two more conditions on the matrix $\Gamma$, namely
\begin{equation}\label{cond-3}
 \pp \left( \|\Gamma\|_{HS} < \sqrt{Nn}/2\r) \leq  p_1,
\end{equation}
for some $p_1\in (0,1)$ and
where $\|\Gamma\|_{HS}$ denotes the Hilbert--Schmidt norm of $\Gamma$; and
\begin{equation}\label{cond-4}
 \pp \left( \|\Gamma\| > \mu \sqrt{N}\r) \leq  p_2,
\end{equation}
for some $p_2\in (0,1)$, $\mu \geq 1$ and where
$\|\Gamma\|$ denotes the operator (spectral) norm of $\Gamma$.
Both conditions are satisfied for example when entries of $\Gamma$
are i.i.d. centered random variables with finite $p$-th moment for some $p>4$. Indeed,
Rosenthal's inequality (see proof of Corollary~6.4 in \cite{GLPT}) implies
(\ref{cond-3}) with $p_1\leq (C_p \E|\xi |^p)/(Nn)^{p/4}$; while Theorem~2.1 combined with
Corollary~6.4 in \cite{GLPT} implies (\ref{cond-4}) with $\mu =C_p'$ and
$$p_2 \leq 1/N^{c_p} + (C_p \E|\xi |^p)N/n^{p/4}$$
(to make $p_2<1$ we have to ask $C_p \E|\xi |^p N\leq  n^{p/4}$).
We would like also to note that the proof below  works also for $N\leq n^{\alpha}$ for some
$\alpha\in (1,2]$  if we substitute
the condition (\ref{cond-4}) with
$$
 \pp \left( \|\Gamma\| > \mu (Nn)^\gamma\r) \leq  p_2
$$
for some $\gamma \in (0, 1/2)$, which could be the case in the absence
of 4-th moment (see for example Corollary~2 in \cite{AAP} and Remark~2
in \cite{LiSp}). Note also that the condition (\ref{cond-4}) implies
(\ref{cond-two}), since $\|\Gamma\| \geq \max_{i\leq N}  |\Gamma ^* e_i| $.

\begin{theorem}
\label{mstar-2}
Let $\mu \geq 1$, $n\geq 16 \mu^2$, and $ 2  n<N\leq n^2$ and assume that the
matrix $\Gamma$ satisfies conditions (\ref{cond-3}) and (\ref{cond-4}) for some
$p_1, p_2\in (0,1)$. Then with probability at least $1-p_1-p_2$
$$
M (K_N^0) \geq c \sqrt{\ln (n/(8 \mu^2))}.
$$

\end{theorem}

\proofs
We apply Vershynin's extension \cite{V} of Bourgain-Tzafriri theorem
\cite{BT}.
Denote $A= \| \Gamma ^* \| _{{HS}}$, $B= \|\Gamma ^* \|$.  Vershynin's theorem
implies that there exists $\sigma
\subset \{1, \ldots, N\}$ of cardinality at least $A^2 /(2 B^2)$
such that for all $i\in \sigma$ one has $| \Gamma ^* e_i | \ge c_3
A/\sqrt{N}$, where $c_3$ is an absolute positive constant, and vectors
$\Gamma ^* e_i$, $i \in \sigma$, are almost orthogonal (up to an
absolute positive constant). Since $\Gamma$ satisfies
conditions (\ref{cond-3}) and (\ref{cond-4}),
with probability at least $1-p_1-p_2$ we have $A\geq \sqrt{Nn}/2$
and $B\leq \mu \sqrt{N}$. Therefore,  there exists
$\sigma \subset \{1, \ldots, n\}$ of cardinality at least
$n/(8 \mu^2)$ such that $| \Gamma ^* e_i | \ge c_3 \sqrt{n}/2$ for $i\in
\sigma$ and $\{ \Gamma ^* e_i \}_{i \in \sigma}$ are almost
orthogonal. Then,
$$
M(K_N^0) \geq \frac{1}{\sqrt{n}} \E \| G \| _{K_N^0}
= \frac{1}{\sqrt{n}} \E \max _{i\leq N} \la G, \Gamma ^*
e_i \ra \geq \frac{1}{\sqrt{n}} \E \max _{i\in \sigma }
\la G, \Gamma ^* e_i \ra .
$$
Since $\{\Gamma ^* e_i\}_{i\in \sigma }$ are almost orthogonal, by
Sudakov inequality (see, e.g., \cite{P}), the last expectation is
greater than $c_4 \sqrt{\ln (n/(8 \mu^2))}$, where $c_4$ is an
absolute constant. This completes the proof.
\kkk

\section{Smallest singular value}
\label{smsing}

In this section we provide a simple short proof of a weaker inclusion,
namely, we obtain a lower bound on the radius of the largest ball inscribed into $K_N$.
It is based on a  lower bound for the smallest singular value for tall
matrices. Although such bounds are known with possibly better
constants (see the last remark in \cite{KM}
or the main theorem of \cite{Tikh}), we would like to emphasize a simple short
proof, based on our Theorem~\ref{Kostya}.
In fact our proof is close to the corresponding proofs in \cite{LPRT} and \cite{LR},
however it is somewhat cleaner and it uses Theorem \ref{Kostya} instead of a
standard net argument via the norm of an operator.
We would also like
to mention that very recently G.~Livshyts has extended such
results to rectangular random matrices with arbitrarily small aspect ratio \cite{Galyna}.


Recall that for an $N \times n$ matrix $\Gamma$ with $N \ge n$,
 its smallest singular value $s_n(\Gamma)$ can be defined by
\[
s_n(\Gamma) = \inf_{x \in S^{n-1}} \|\Gamma x \|_2.
\]
In this section we assume that the random matrix $\Gamma$ satisfies
conditions described at the beginning of Section~\ref{geom} with fixed
$u, v\in (0,1)$. Recall that $c_{uv}=c u v \sqrt{1- v}$ is the constant
from Lemma~\ref{levy}. It will be also convenient to fix two more constants
depending on $v$,
\begin{equation*}
 \gam_1 = \gam_1 (v)  :=
  \left\{
\begin{array}{ll}
 \sqrt{\ln 2} & \mbox{ if }\,\, v\geq 1/2,
\\
 \sqrt{\ln \frac{1}{v}}  & \mbox{ if }\,\,  v<1/2
\end{array}
\right.
\quad \mbox{ and } \quad
\gam_2 = \gam_2 (v)  :=
  \left\{
\begin{array}{ll}
  \ln\frac{2}{1+v} & \mbox{ if }\,\, v\geq 1/2,
\\
 \ln \frac{1}{2v-v^2}  & \mbox{ if }\,\,  v<1/2.
\end{array}
\right.
\end{equation*}

\begin{theorem}\label{sing-val}
There exist an absolute  constant $C_0>1$
such that for  $N\geq \big(\frac{C_0}{\gam_2}\, \ln \frac{1}{c_{uv}}\big)\, n$ one has
$$
  \pp \left( s_n (\Gamma) \le \frac{c_{uv}\sqrt{\gam_2}}{4\gam_1} \, \sqrt{N}  \r) \leq
  3 \exp \left( - \min\{2, \gam_2\} N/8 \r).
$$
\end{theorem}

\bigskip

Since
$$
 h_{\Gamma ^* B_1^N} (x) =\| \Gamma ^* x\|_\infty
 \quad \mbox{ and } \quad
   K_N = \Gamma ^* B_1^N \supset \frac{1}{\sqrt N} \Gamma ^* B_2^N,
$$
 this theorem immediately implies the following inclusion.

\begin{sled} \label{incl-ball}
For  $N\geq \big(\frac{C_0}{\gam_2}\, \ln \frac{1}{c_{uv}}\big)\, n$ one has
$$
  \pp \left( K_N \supset \frac{c_{uv}\sqrt{\gam_2}}{4\gam_1}  \sqrt{N} B_2^n\r)
   \geq 1- 3 \exp \left( - \min\{2, \gam_2\} N/8 \r).
$$
\end{sled}

\bigskip



To prove Theorem~\ref{sing-val} we first provide the individual bounds.

\begin{prop}
\label{propo-tall}
Let $1 \leq n < N$.
 Then for every $x \in S^{n-1}$  one has
$$
\P \Bigl( \|\Gamma x\|_2 \leq \frac{c_{uv}\sqrt\gam_2}{2\gam_1}\, \sqrt{N}\Bigr)
\leq  \exp \left(- 3 \gam _2 N/4 \r).
$$
\end{prop}

\bigskip

%
%

\proofs
Fix  $x = (x_1, \ldots, x_n) \in \R$ with $\euclidnorm{x} = 1$.
Denote $f_j := |\sum_{i=1}^{n} \xi_{ji} x_i|$, so that
$$
  \euclidnorm{\Gamma x}^2 = \sum_{j=1}^{N} f_j^2.
$$
Clearly $f_1,\ldots,f_N$ are independent. Therefore, for any $t, \tau > 0$ one has
\begin{align*}
\P \bigl( \euclidnorm{\Gamma x}^2 \leq t^2 N \bigr)
&= \P \biggl( \sum_{j=1}^{N} f_j^2 \leq t^2 N \biggr)
= \P \biggl( \tau N - \frac{\tau}{t^2} \sum_{j=1}^{N} f_j^2 \geq 0 \biggr) \notag\\
&\leq \EE \exp \biggl( \tau N - \frac{\tau}{t^2} \sum_{j=1}^{N} f_j^2 \biggr)
= e^{\tau N} \prod_{j=1}^{N} \EE \exp \biggl( -\frac{\tau f_j^2}{t^2} \biggr).
\end{align*}
Lemma~\ref{levy} implies that $\P (f_j < c_{uv}) \leq v$ for every $j\leq N$.
Write $\tau =t^2\eta/c_{uv}^2$ for some $\eta >0$. Then
\begin{align*}
\EE \exp \Bigl( -\frac{\tau f_j^2}{t^2} \Bigr)
&= \int_{0}^{1} \P \biggl( \exp \Bigl(  -\frac{\eta f_j^2}{c_{uv}^2} \Bigr) > s \biggr)\, ds\\
&= \int_{0}^{e^{-\eta}} \!\!\! \P \biggl( \exp \Bigl(  \frac{\eta f_j^2}{c_{uv}^2} \Bigr) < \frac{1}{s} \biggr)\, ds +
\int_{e^{-\eta}}^{1} \!\!\! \P \biggl( \exp \Bigl( \frac{\eta f_j^2}{c_{uv}^2} \Bigr) < \frac{1}{s} \biggr)\, ds\\
&\leq e^{-\eta} + \P (f_j < c_{uv}) ( 1 - e^{-\eta} )
\leq
e^{-\eta} + v ( 1 - e^{-\eta} )  .
\end{align*}
Choose $\eta =\gam _1^2 =\ln \max\{2, 1/v\}$. Then the right hand side is $e^{-\gam_2}$.
Therefore
\begin{equation*}
  \P \bigl( \euclidnorm{\Gamma x}^2 \leq t^2 N \bigr)
  \leq e^{\tau N} e^{-\gam_2 N }
  = \exp \left(-N (\gam _2 -  t^2\gam _1^2/c_{uv}^2) \r).
\end{equation*}
Choosing $t=\sqrt{\gam_2}  c_{uv}/(2\gam _1)$ we complete the proof.
\qed

\bigskip

\noindent
{\bf Proof of Theorem~\ref{sing-val}. }
Let $\delta =\min\{1, \gam_2/2\}$. Note
that  $n/(2N)\leq \delta \leq 1$.
Let $C\geq 1$ be the absolute constant from  Theorem~\ref{Kostya}.
Set
$$
 \eps := \frac{c_{uv}\,  \sqrt{\gam _2 \delta }}{4C \gam _1 \sqrt{ n }}  <\frac{1}{ \sqrt{ n }}.
$$
By Theorem~\ref{Kostya} (see Remark~\ref{rem:appli}), applied with $T = S^{n-1}$ and $k = N$,
there  exists a net $\mathcal{N}\subset B_2^n$  with cardinality at most
$$
  \left(\frac{224\delta  N}{\eps n^{3/2}}\r)^n\,  e^{\delta N}\leq
  \left(\frac{896 C\gam _1 \sqrt{\delta} N}{c_{uv} \sqrt{\gam_2} n} \r)^n\,  e^{\delta N}
$$
such that with probability at least $1- e^{-\delta N/4} - e^{-N}$ one has
$$
 \forall x\in B_2^n \,\, \exists y_x\in \mathcal{N}
  \, \, \, \, \, \, \, \, \mbox{ such that } \, \, \, \,\, \, \, \,
  \Gamma (x-  y_x) \in  C\eps  \sqrt{ N n/\delta } \, B_2^n =(c_{uv}\sqrt{\gam_2}/(4\gam_1)) \sqrt{N} \, B_2^n.
$$
Condition on the corresponding event, denoted below by $\Omega_{0}$. Assume that
$x \in S^{n-1}$ satisfies $\euclidnorm{\Gamma x} \leq (c_{uv}\sqrt{\gam_2}/(4\gam_1))  \sqrt{N}$. Then
for the corresponding $y_x\in \mathcal{N}$ we have
\begin{equation*}
 \euclidnorm{\Gamma y_x} \leq \euclidnorm{\Gamma x} + \euclidnorm{\Gamma(y_x - x)}
 \leq
 (c_{uv}\sqrt{\gam_2}/(2\gam_1))  \sqrt{N}.
\end{equation*}
This implies
$$
 q_0:=  \P \left(\exists x\in S^{n-1} \, \, | \,\, \euclidnorm{\Gamma x} \leq \frac{c_{uv}\sqrt{\gam_2}}{4\gam_1}  \sqrt{N}\r)
  \leq \P \left(\Omega_{0}^c \r)
  +\P \left(\exists y\in \mathcal{N} \, \, | \,\, \euclidnorm{\Gamma y} \leq \frac{c_{uv}\sqrt{\gam_2}}{2\gam_1}  \sqrt{N} \r) .
$$
Applying Proposition~\ref{propo-tall} and using $\delta \leq \gam_2/2$,
$$
  q_0 \leq 2 e^{-\delta N/4} +   \left(\frac{896 C\gam _1 \sqrt{\delta} N}{c_{uv} \sqrt{\gam_2} n} \r)^n\,
  \exp \left(-\gam _2 N/4  \r).
$$
Using formulas for $c_{uv}$, $\gam_1$, $\gam_2$, and $\delta$,
it is not difficult to check that there exists
an absolute constant $C_1>0$ such that
$$
  \ln \frac{896 C\gam _1 \sqrt{\delta} }{c_{uv} \sqrt{\gam_2} }
  \leq C_1 \ln \frac{1}{c_{uv}} .
$$
Therefore there exists another absolute constant $C_2>0$ such that
$$
 \left(\frac{896 C\gam _1 \sqrt{\delta} N}{c_{uv} \sqrt{\gam_2} n} \r)^n\,
  \exp \left(-\gam _2 N/4  \r)\leq \exp \left(-\gam _2 N/8  \r),
$$
provided that
$$
  N/n\geq (C_2/\gam_2)\, \ln (1/c_{uv}) .
$$
This completes the proof.
\qed


\address




\begin{thebibliography}{99}



\bibitem{AAP}
A. Auffinger, G. Ben Arous, S. P\'ech\'e,
{\em  Poisson convergence for the largest
eigenvalues of heavy tailed random matrices,}
Ann. Inst. Henri Poincar Probab. Stat. 45
(2009), 589--610.






\bibitem{BF} I.~B\'ar\'any, Z.~F\"uredy, {\em Approximation of the
  sphere by polytopes having few vertices,} Proc. Amer. Math. Soc.
102 (1988), no. 3, 651--659.

\bibitem{Barany} I.~B\'ar\'any and A.~P\'or, {\em On 0-1 Polytopes
 with many facets,} Adv. Math. 161 (2001), 209--228.

\bibitem{BM} J.~Bourgain and  V.~D.~Milman, {\em New volume ratio
  properties for symmetric bodies in $\R$}, Invent. Math. 88 (1987),
no 2, 319--340.
%
\bibitem{BT} J. Bourgain and L. Tzafriri, {\em Invertibility of
  "large" submatrices with applications to the geometry of Banach
  spaces and harmonic analysis}, Israel J. Math. 57 (1987),
137--224.
%
%
\bibitem{CP} B.~Carl and A.~Pajor, {\it Gelfand numbers of
  operators with values in a Hilbert space,} Invent. Math. {\bf 94}
(1988), 479--504.
%


%
%
%



\bibitem{Furedi} M.E.~Dyer, Z.~F\"uredi and C.~McDiarmid, {\em
 Volumes spanned by random points in the hypercube,} Random
Structures Algorithms 3 (1992), 91--106.


\bibitem{GH} A.~Giannopoulos, M.~Hartzoulaki, {\em Random spaces
  generated by vertices of the cube}, Discrete Comp. Geom., 28
(2002), 255--273.



%
\bibitem{gldistance} E.D. Gluskin, {\em The diameter of Minkowski
  compactum roughly equals to $n$},  Funct. Anal. Appl.,~15
(1981), 57--58 (English translation).
%
\bibitem{Gl} E.D. Gluskin, {\em Extremal properties of orthogonal
  parallelepipeds and their applications to the geometry of Banach
  spaces}, (Russian) Mat. Sb. (N.S.)  136 (178) (1988), no. 1,
85--96; translation in Math. USSR-Sb. 64 (1989), no. 1, 85--96.
%
%

\bibitem{Gl3} E.D. Gluskin,
%
{\em The octahedron is badly approximated by random subspaces},
Funct. Anal. Appl. 20 (1986), 11--16; translation from Funkts. Anal. Prilozh.
20 (1986), no. 1, 14--20.

\bibitem{GLSW}
Y. Gordon, A.E. Litvak, C. Schuett, E. Werner,
{\em Geometry of spaces between zonoids and polytopes,} Bull. Sci. Math., 126 (2002), 733--762.



\bibitem{GGMP}
Y. Gordon, O. Gu\'edon, M. Meyer, A. Pajor,
{\em
 Random Euclidean sections of some classical Banach spaces},
 Math. Scand. 91 (2002), no. 2, 247--268.


\bibitem{OG}
O. Gu\'edon, {\em
Gaussian version of a theorem of Milman and Schechtman},
Positivity 1 (1997), no. 1, 1--5.

\bibitem{GLPT}
O. Gu\'edon, A.E. Litvak, A. Pajor, N. Tomczak-Jaegermann,
{\em On the interval of fluctuation of the singular values of random matrices,}
J. Eur. Math. Soc., 19 (2017), 1469--1505.




%
%
%

\bibitem{KM}
V. Koltchinskii, S. Mendelson,
{\em  Bounding the smallest singular value of a random
matrix without concentration,}
Int. Math. Res. Not. 2015, No 23,   12991--13008.






\bibitem{KKR} F. Krahmer, C. Kummerle, and H. Rauhut,
{\em A Quotient Property for Matrices with Heavy-Tailed Entries and its
Application to Noise-Blind Compressed Sensing}, preprint,
arXiv:1806.04261.







\bibitem{LPRT} A.E. Litvak, A. Pajor, M. Rudelson,
N.~Tomczak-Jaegermann, {\em
Smallest singular value of random matrices and geometry of random polytopes,}
Adv. Math., 195 (2005), 491--523.




\bibitem{LR} A.E. Litvak, O. Rivasplata,
{\em Smallest singular value of sparse random matrices,}
Stud. Math., 212 (2012), 195--218.



\bibitem{LiSp}
A.E. Litvak, S. Spektor,
{\em Quantitative version of a Silverstein's result,}
GAFA, Lecture Notes in Math., 2116 (2014), 335--340.


\bibitem{Galyna}
G.V. Livshyts,
{\em The smallest singular value of heavy-tailed not necessarily
i.i.d. random matrices via random rounding,}
arXiv: 1811.07038.




%
\bibitem{mat} P. Mankiewicz and N. Tomczak-Jaegermann, {\em
  Quotients of finite-dimensional Banach spaces; random phenomena.}
In:``Handbook in Banach Spaces'' Vol II, ed. W.~B.~Johnson,
J.~Lindenstrauss, Amsterdam: Elsevier (2003), 1201--1246.
%
%
%
\bibitem{MS} S.J. Montgomery-Smith, {\em The distribution of
    Rademacher sums}, Proc. Amer. Math. Soc. 109 (1990), no. 2,
  517--522.
%
\bibitem{P} G. Pisier, {\em The Volume of Convex Bodies and Banach
  Space Geometry}, Cambridge University Press, Cambridge, 1989.
%

\bibitem{RT}  E. Rebrova, K. Tikhomirov,
{\em Coverings of random ellipsoids, and invertibility of matrices with i.i.d.
heavy-tailed entries}, Isr. J. Math., {\bf 227}  (2018), 507--544.



\bibitem{ROG}
B.A. Rogozin,
{\em On the increase of dispersion of sums of independent random variables},
Teor. Verojatnost. i Primenen {\bf 6} (1961), 106--108.





\bibitem{CS}
C. Sch\"utt,
{\em
Entropy numbers of diagonal operators between symmetric Banach spaces},
J. Approximation Theory 40 (1984), 121--128.

\bibitem{Tikh}
 K.E. Tikhomirov,
 {\em The smallest singular value of random rectangular matrices with
 no moment assumptions on entries}, Isr. J. Math. 212 (2016),  289--314.


\bibitem{silv} J.W. Silverstein, {\em The smallest eigenvalue of
  a large dimensional Wishart matrix}, Ann. Probab. 13 (1985),
1364--1368.
%
%
\bibitem{szbasis} S.J.~Szarek, {\em The finite-dimensional basis
  problem with an appendix on nets of Grassman manifold}, Acta Math.
141 (1983), 153--179.
%
%
\bibitem{V} R.~Vershynin, {\em John's decompositions: selecting a
  large part}, Israel J. Math. 122 (2001), 253--277.
%
%
%
\bibitem{Ziegler} G.~M.~Ziegler, {\em Lectures on 0/1 polytopes,}
in ``Polytopes-Combinatorics and Computation" (G. Kalai and G. M.
Ziegler, Eds), pp. 1--44, DMV Seminars, Birkh\"auser, Basel, 2000.
%
\end{thebibliography}
\end{document}